\documentclass[10pt]{article}
\usepackage{graphicx}
\usepackage{amssymb}
\usepackage{epstopdf}
\usepackage{amsmath}
\usepackage{amsfonts}

\def\phi{{\varphi}}

\DeclareSymbolFont{AMSb}{U}{msb}{m}{n}
\DeclareMathSymbol{\N}{\mathbin}{AMSb}{"4E}
\DeclareMathSymbol{\Z}{\mathbin}{AMSb}{"5A}
\DeclareMathSymbol{\R}{\mathbin}{AMSb}{"52}
\DeclareMathSymbol{\Q}{\mathbin}{AMSb}{"51}
\DeclareMathSymbol{\I}{\mathbin}{AMSb}{"49}
\DeclareMathSymbol{\C}{\mathbin}{AMSb}{"43}

\def\be{\begin{equation}}
\def\ee{\end{equation}}
\def\ber{\begin{eqnarray}}
\def\eer{\end{eqnarray}}

\def\beq{\begin{equation}}
\def\eeq{\end{equation}}

\begin{document}

\addtolength{\textheight}{0 cm} \addtolength{\hoffset}{0 cm}
\addtolength{\textwidth}{0 cm} \addtolength{\voffset}{0 cm}

\newenvironment{acknowledgement}{\noindent\textbf{Acknowledgement.}\em}{}

\setcounter{secnumdepth}{5}
 \newtheorem{proposition}{Proposition}[section]
\newtheorem{theorem}{Theorem}[section]
\newtheorem{lemma}[theorem]{Lemma}
\newtheorem{coro}[theorem]{Corollary}
\newtheorem{remark}[theorem]{Remark}
\newtheorem{claim}[theorem]{Claim}
\newtheorem{conj}[theorem]{Conjecture}
\newtheorem{definition}[theorem]{Definition}
\newtheorem{application}{Application}

\newtheorem{corollary}[theorem]{Corollary}

\title{  Multi-marginal Monge-Kantorovich  transport problems: A characterization of   solutions  \\
 Probl\`emes de transport multi-marginal de Monge-Kantorovich: Une caract\'erisation des solutions}
\author{
Abbas Moameni
\thanks{Supported by a grant from the Natural Sciences and Engineering Research Council of Canada.}
\hspace{2mm}\\
{\it\small School of Mathematics and Statistics}\\
{\it\small Carleton University}\\
{\it\small Ottawa, ON, Canada K1S 5B6}\\
{\it\small  momeni@math.carleton.ca}\\
%
}

\maketitle

\begin{abstract}
We shall present  a measure theoretical approach   for which together with  the Kantorovich duality provide an efficient    tool to study the
optimal  transport problem. Specifically,
 we study the support of optimal plans
 where the  cost function does not satisfy the classical twist condition in the two marginal problem as well as in  the multi-marginal 
case when twistedness is limited to certain subsets. 
\end{abstract}
\section*{R\'esum\'e}
Dans cet article, nous \'etudions le probl\`eme de transport optimal du
point de vue de la th\'eorie de la mesure, \`a l'aide de la  dualit\'e de
Kantorovich. En particulier, nous \'etudions le support des plans optimaux o\`u la
fonction co\^ut ne satisfait pas la condition
de ``twist'' dans le probl\`eme \`a deux marginales, 
 ainsi que dans le cas multi-marginale quand la condition 
``twist''  est limit\'ee \`a des sous-ensembles pr\'ecis.\\

\section{Introduction}
We consider the Monge-Kantorovich transport problem  for Borel probability measures $\mu_1, \mu_2,..., \mu_n$ on smooth manifolds    $X_1, X_2,..., X_n.$
The cost function $c: X_1 \times X_2 \times... \times X_n  \to \R$ is bounded and continuous. Let  $\Pi(\mu_1,...,\mu_n)$ be  the set of Borel probability measures on
$X_1 \times X_2 \times... \times X_n$  which have $X_i$-marginal $\mu_i$ for each $i \in \{1,2,...,n\}.$
The transport cost associated to a transport plan
 $\pi \in \Pi(\mu_1,...,\mu_n)$ is given by
\[I_c(\pi)=\int_{X_1 \times X_2 \times... \times X_n} c(x_1,...,x_n) \, d \pi.\]
We consider  the Monge-Kantorovich transport problem,
\[  \inf \{ I_c(\pi); \pi \in \Pi(\mu_1,...,\mu_n)\}. \qquad \qquad (MK)        \]
If  a transport plan minimizes the cost, it will be called an optimal plan.
We say that an  optimal plan  $\gamma$ induces a Monge solution if it is concentrated on the graph $\{(x, T(x); \, x \in X_1\}$ of a measurable map 
$T : X_1  \to X_2 \times... \times X_n.$ 
 In contrary to the Monge problem, the Kantorovich problem always admits solutions as soon as the
cost function is a non-negative lower semi continuous  function (see \cite{V} for a proof). 
When $n=2$, a  general criterion for existence and uniqueness of an optimal transport map known as the twist condition
  dictates  the map
$  y \to D_1 c(x,y) $
to be  injective for fixed $x \in X_1.$  Under the twist condition and the absolute continuity of $\mu_1$, the  optimal plan $\gamma$    which solves the
Monge-Kantorovich  problem $(MK)$ is supported on the graph of an optimal transport map $T,$ i.e., $\gamma=(\text{Id} \times T)_\#\mu.$ For larger $n$, questions 
regarding the existence and uniqueness 
are not fully understood yet. By now there are  many interesting results for the multi-marginal problem in the general case as well as particular models
 (see for instance \cite{B-D-G, Ca, C-M-N, G-S, G-M, He, K-P, P1,P2},  the bibliography
is not exhaustive).
 When $n>2,$ as shown in \cite{K-P}, the twist 
condition can be replaced by twistedness on $c$-splitting sets.
\begin{definition} A set $S \subset  X_1 \times X_2\times ...\times X_n$ is a $c$-splitting set if there exists
Borel functions $u_i : X_i \to \R$  such that for all  $(x_1, x_2, ..., x_n),$
\[ \sum_{i=1}^n  u_i(x_i) \leq  c(x_1, x_2, ..., x_n) \]
with equality whenever $(x_1, x_2, ..., x_n) \in S$.  The $n$-tuple $(u_1,...,u_n)$ is called the $c$-splitting tuple for $S.$
\end{definition}
 In \cite{Moa}, for the case $n=2,$ the author  relaxed the twist condition by a new property, i.e., 
\begin{itemize}
\item {{Generalized-twist condition:}}   We say that $c$  satisfies the generalized-twist condition if
for any $\bar x_1 \in X_1$ and  $\bar x_2 \in X_2$
  the set $L_{(\bar x_1,\bar x_2)}:=\{x_2 \in X_2  ;\, D_1 c(\bar x_1,x_2)=D_1 c(\bar x_1, \bar x_2)\}$ 
is a finite subset of $X_2$. Moreover,  if there exists $m \in \mathbb{N}$ such that  for each $\bar x_1 \in X_1$ and $\bar x_2 \in X_2$ the cardinality of the set $L_{(\bar x_1,\bar x_2)}$  does 
not exceed $m$ then we say that $c$ satisfies the 
$m$-twist condition. 
\end{itemize}
 Under the $m$-twist condition, it is shown that for  each optimal plan
  $\gamma $ of  $(MK),$
 there exist a sequence of non-negative measurable real functions 
$\{\alpha_i\}_{i=1}^m $  on  $X_1$ with $\sum_{i=1}^m \alpha_i =1$
 and,  Borel measurable maps $G_1,...,G_m: X_1 \to  X_2$ such that
$
\gamma =\sum_{i=1}^m \alpha_i (\text{Id} \times  G_i)_\# \mu.
$\\
Our aim  in this work is to extend this result to the multi-marginal case.  For the rest of the  paper we always assume that $c$ is  non-negative,  lower semi-continuous,  $\oplus_{i=1}^n \mu_i$-a.e. differentiable with respect to the first variable and that $I_c(\gamma)$ is finite for some transport plan $\gamma.$ We also assume that the Kantorovich dual problem admits a solution 
$(\phi_1,...,\phi_n)$   such that $\phi_1$ is differentiable $\mu_1$-a.e., 
 $\phi_1(x_1)+...+\phi_n(x_n) \leq c(x_1,..., x_n)$ for all $(x_1,...,x_n)$ and
 \begin{equation*}
 \int c \, d\gamma= \sum_{i=1}^n\int_{X_i} \phi_i(x_i) \, d \mu_i.
 \end{equation*}

 We denote by $D_1(c)$ the set of points at which $c$ is  differentiable with respect to the first variable.  The generalized twist structure takes the following form in the multi-marginal case.
\begin{definition} Let  $c$ be a Borel measurable function. \\
1.  {{\bf $\mathbf{m}$-twist condition:}} 
Say that 
$c$ is $m$-twisted on $c$-splitting sets 
if  for any $c$-splitting set  $S \subset  X_1 \times X_2 \times  ...\times X_n$ and  any $(\bar x_1, \bar x_2, ..., \bar x_n)\in S\cap D_1(c)$ 
 the cardinality of the set 
\[\Big\{(\bar x_1,x_2, ..., x_n) \in S\cap D_1(c);\,D_{x_1}c(\bar x_1, \bar x_2, ..., \bar x_n)= D_{x_1}c(\bar x_1,  x_2, ...,  x_n) \Big\},\]
is at most $m.$\\
2. {{\bf Generalized-twist condition:}} Say that 
$c$ satisfies the generalized twist condition on $c$-splitting sets if for any
 $c$-splitting set  $S \subset  X_1 \times X_2 \times  ...\times X_n$ and  any $(\bar x_1, \bar x_2, ..., \bar x_n)\in S\cap D_1(c)$ 
  the set 
\[\Big\{(\bar x_1,x_2, ..., x_n) \in S\cap D_1(c);\,D_{x_1}c(\bar x_1, \bar x_2, ..., \bar x_n)= D_{x_1}c(\bar x_1,  x_2, ...,  x_n) \Big\},\]
is a finite subset of $S$.
\end{definition}
The following result provides a connection between the generalized twist condition and the local $1$-twistedness.
\begin{proposition}\label{non}
Assume that $c$ is continuously differentiable with respect to the first variable and $S$ is a compact  $c$-splitting set. 
 If $c$ is locally $1$-twisted on $S$   then $c$ satisfies the generalized-twist condition on $S$.  
\end{proposition}

We now state our main result in this paper. 
\begin{theorem} \label{main} Assume that  the cost function $c$   satisfies  the $m$-twist condition on $c$-splitting sets and 
$\mu_1$ is non-atomic. 
  Then for each optimal plan
  $\gamma $ of  $(MK)$ with $Supp(\gamma) \subset D_1(c),$
 there exist $k \leq m,$ a sequence of non-negative measurable real functions 
$\{\alpha_i\}_{i=1}^k $ on $X_1$ 
 and,  Borel measurable maps $G_1,...,G_k: X_1 \to  X_2 \times ...\times X_n$ such that
\begin{equation}\label{late}\gamma =\sum_{i=1}^k \alpha_i (\text{Id} \times  G_i)_\# \mu, \qquad \end{equation}
where  $\sum_{i=1}^k \alpha_i(x) =1$ for $\mu_1$-a.e. $x \in X_1$
\end{theorem}

 We also have the following for the generalized-twist condition.
\begin{proposition}\label{mm}
 Under the assumption of Theorem \ref{main}, if one replaces the $m$-twist condition by the generalized -twist condition then   each optimal plan
  $\gamma $ of  $(MK)$ is of the form (\ref{late}) with $k \in \mathbb{N}\cup \{\infty \}.$
\end{proposition}
As shown in \cite{Moa}, the most interesting examples of costs satisfying the generalized-twist condition are non-degenerate costs on smooth $n$-dimensional manifolds $X$ and $Y$. 
  Denote by $D^2_{xy} c(x_0,y_0)$
the $n\times n$ matrix of mixed second order partial derivatives of the function $c$ at the point $(x_0,y_0).$
A cost $c \in C^2(X \times Y)$ is non-degenerate provided
$D^2_{xy} c(x_0,y_0)$ is non-singular, that is $\text{det} \Big( D^2_{xy} c(x_0,y_0)\Big)\not=0.$ In our forthcoming work \cite{M-P},
 following an idea in \cite{P2} together with Proposition \ref{non}, a 
differential condition similar to  the non-degeneracy condition (in $n=2$) is derived for the multi-marginal case that guaranties the general 
twist property on $c$-splitting sets and consequently the characterization of the support of optimal plans due to Theorem \ref{main}.\\
In the next section, we shall discuss the key ingredients  for our methodology in this work. Section 3 is devoted to  the proof  of the main results,
\section{Measurable  sections and extremality}
 Let $(X,  \mathcal{B}, \mu)$ be a finite, not necessarily complete measure space, and $(Y, \Sigma)$  a measurable space.
 The completion of $\mathcal{B}$
with respect to $\mu$ is denoted by $\mathcal{B}_\mu,$ when necessary, we identify $\mu$ with its
completion on $\mathcal{B}_\mu.$ 
 The push forward of the measure $\mu$ by a  map $T: (X,  \mathcal{B}, \mu) \to (Y, \Sigma)$ is denoted by
$T_\# \mu,$ i.e.
\[T_\# \mu (A)=\mu(T^{-1}(A)), \qquad \forall A \in \Sigma.\]

\begin{definition} Let $T: X \to Y$  be $(\mathcal{B}, \Sigma)$-measurable and $\nu$ a positive measure on $\Sigma.$
 We call a map  $F: Y \to X$  a $(\Sigma_\nu,\mathcal{B})$-measurable  section of $T$ if $F$ is $(\Sigma_\nu,\mathcal{B})$- measurable
 and  $T \circ F = \text{Id}_Y.$
\end{definition}

If $X$ is a topological space we denote by $\mathcal{B}(X)$ the set of Borel sets on
$X.$ The space of Borel probability measures on a topological space $X$ is denoted by $\mathcal{P}(X)$.
   For a measurable map $T: (X, \mathcal{B}(X) ) \to (Y,\Sigma, \nu)$ denote by $\mathcal{M}( T, \nu)$  
the set of all  measures $\lambda$ on $\mathcal{B} $ so that $T$ pushes $\lambda$ forward to $\nu,$ i.e.
\[\mathcal{M}( T, \nu)=\{\lambda \in \mathcal{P}(X) ; \, T_\# \lambda=\nu\}.\]
Evidently $\mathcal{M}( T, \nu)$ is a convex set. A measure $\lambda$ is an extreme point of $\mathcal{M}( T, \nu)$ if  the identity $ \lambda= \theta \lambda_1+(1-\theta )\lambda_2$
with $\theta \in (0,1) $ and $\lambda_1, \lambda_2 \in
\mathcal{M}( T, \nu)$ imply that  $ \lambda_1=\lambda_2$.
 The set of extreme points of  $\mathcal{M}( T, \nu)$ is denoted by $ext\,\mathcal{M}( T, \nu).$\\

    We recall the following result from \cite{Gr} in
 which a characterization of the set $ext\,\mathcal{M}( T, \nu)$ is given.
\begin{theorem}\label{Graf} Let $(Y,\Sigma, \nu)$ be a probability space, $(X, \mathcal{B}(X))$ be a Hausdorff space with a Radon probability
measure $\lambda$,
and let $T : X\to Y$ be an $(\mathcal{B}(X), \Sigma)$-measurable mapping.
Assume that  $T$ is surjective and $\Sigma$ is countably separated. The following conditions are equivalent:\\
(i) $\lambda$ is an extreme point of $M(T, \nu)$;\\
(ii) there exists a $(\Sigma_\nu,\mathcal{B}(X))$-measurable section $F : Y \to X$ of the mapping
$T$ with $\lambda = F_\# \nu$.\\
\end{theorem}
By making  use of the Choquet theory  in the setting of noncompact sets of measures \cite{WW}, each  $ \lambda \in M(T, \nu)$ can be represented  as a Choquet type integral over
 $ ext \, M(T, \nu).$
Denote by $\Sigma_{ ext \, M(T, \nu)}$ the $\sigma$-algebra over  $ext\, M(T, \nu)$ generated by  the functions $\varrho \to \varrho(B),$  $B \in  \mathcal{B}(X).$
We have the following result (see \cite{Moa} for a proof).

\begin{theorem} \label{cho}
Let $X$ and $Y$ be complete separable metric spaces and $\nu$ a probability measure on $\mathcal{B} (Y).$
Let $T:(X, \mathcal{B}(X)) \to (Y, \mathcal{B}(Y))$ be a surjective measurable mapping and  let $\lambda \in M(T, \nu).$ Then there exists a  probability measure $\xi$ on $\sum_{ext \,  M(T, \nu)}$
 such that for each $B \in \mathcal{B}(X)$,

\[\lambda(B)=\int_{ext \, M(T, \nu)} \varrho(B)\, d\xi(\varrho), \qquad  \big (\varrho \to \varrho(B) \text{ is measurable}\big).\]
\end{theorem}
 \section{Proofs.}
 In this section we shall proceed with the proofs of  the statements in the introduction. We first state  some  preliminaries required for the proofs. Let $ \gamma$ be a solution of $(MK)$ such that $Supp(\gamma) \subset D_1(c).$
  It is standard that $\gamma \in \Pi (\mu_1,..., \mu_n)$ is non-atomic if and only if  at least one  $\mu_i$ is non-atomic (see for instance \cite{R}).
 Set $Y= X_2 \times ...\times X_n.$
Since $\mu_1$ is non-atomic
it follows  that the Borel measurable spaces  $(X_1, \mathcal{B}(X_1), \mu_1)$  and $(X_1 \times  Y, \mathcal{B}(X_1 \times Y), \gamma)$ are  isomorphic.
Thus, there exists an isomorphism $T=(T_1, T_2, ..., T_n)$ from $(X_1, \mathcal{B}(X_1), \mu_1)$  onto $(X_1 \times Y, \mathcal{B}(X_1 \times Y), \gamma)$.
 It can be easily deduced that
  $T_i: X_1 \to X_i$ are surjective maps and
\[T_i \# \mu_1=\mu_i, \qquad i=1,2,...,n.\]
Consider the convex set \[\mathcal{M}(T_1, \mu_1)=\big\{ \lambda \in \mathcal{P}(X_1); \, T_1\# \lambda=\mu_1\big \},\] and note that $\mu_1 \in \mathcal{M}(T_1, \mu_1).$
The following definition and proposition are  essential in the sequel (see  \cite{Moa} for a proof).
\begin{definition}\label{gen} Denote by $\mathcal{S}(T_1)$ the set of all sections of $T_1$.
Let $\mathcal{K} \subset \mathcal{S}(T_1).$ We say that a  measurable function  $F:\big (X, \mathcal{B}(X)_\mu\big ) \to \big (X,\mathcal{B}(X)\big )$ is generated by $\mathcal{K}$ if there exist a sequence $\{F_i\}_{i=1}^\infty \subset \mathcal{K}$ such that $X=\cup_{i=1}^\infty A_i$ where 
\[A_i= \{x \in X; \, \, F(x)=F_i(x)\}.\]
We also denote by $\mathcal{G}(\mathcal{K})$ the set of all functions generated by $\mathcal{K}.$  It is easily seen that $\mathcal{K} \subset \mathcal{G}(\mathcal{K})\subset \mathcal{S}(T_1).$ 
\end{definition}

\begin{proposition} \label{vint} Let $\mathcal{K}$ be a nonempty subset of $\mathcal{S}(T_1).$ 
Then there exist $ k \in \mathbb{N} \cup\{\infty\}$ and a   sequence $\{F_i\}_{i=1}^k \subset \mathcal{G}(\mathcal{K})$  such that  the following assertions hold:
\begin{itemize}
\item[i.] for each $i\in \mathbb{N}$ with $i\leq k$ we have $\mu(B_{i})>0$  where $\{B_i\}_{i=1}^k$ is defined recursively as follows
\[B_1=X \quad \& \quad B_{i+1}=\Big \{x \in B_i; \,\, F_{i+1}(x)\not\in \{F_1(x),..., F_i(x) \}\Big\}\quad \text{provided }  k>1.\]
\item[ii.]  For all $F \in \mathcal{G}(\mathcal{K})$  we  have 
\[\mu \Big  \{x \in B_{i+1}^c\setminus B^c_{i}; \, \, F(x) \not\in \{F_1(x),..., F_i(x) \}\Big\}=0.\]
\item[iii.]If $k\not =\infty$ then for all $F \in \mathcal{G}(\mathcal{K})$ 
\[\mu \Big   \{x \in B_k; \, \, F(x) \not\in \{F_1(x),..., F_k(x) \}\Big\}=0.\]
\end{itemize}
Moreover, if either $k\not=\infty$ or $k=\infty$ and $\mu(\cap_{i=1}^\infty B_i)=0$ then for every $F \in \mathcal{G}(\mathcal{K})$ the measure $\varrho_F=F_\#\mu$ is absolutely continuous with respect to the measure $\sum_{i=1}^k \varrho_i $ where $\varrho_i={F_i}_\#\mu.$ 
\end{proposition}

\textbf{Proof of Theorem \ref{main}.}
Since $ \mu_1 \in  \mathcal{M}(T_1, \mu_1)$,  it follows from  Theorem \ref{cho}  that   there exists a  probability measure
$\xi$ on $\sum_{ext \,  M(T_1, \mu_1)}$
 such that for each $B \in \mathcal{B}(X_1)$,

\begin{equation}\label{kre} \mu(B)=\int_{ext \, M(T_1, \mu_1)} \varrho(B)\, d\xi(\varrho), \qquad  \big (\varrho \to \varrho(B) \text{ is measurable}\big).\end{equation}

On the other hand, there exist   functions $\{\phi_i\}_{i=1}^n$   such that
 $\phi_1(x_1)+...+\phi_n(x_n) \leq c(x_1,..., x_n),$ $\phi_1$ is $\mu_1$-a.e. differentiable, and that   
 \begin{equation*}
 \int c \, d\gamma= \sum_{i=1}^n\int_{X_i} \phi_i(x_i) \, d \mu_i.
 \end{equation*}
Let $S$ be the $c$-splitting set generated by the $n$-tuple $(\phi_1,...,\phi_n)$, that is, \[S=\{(x_1,...,x_n); \, c(x_1,..., x_n)=\sum_{i=1}^n\phi_i(x_i)  \}.
\] As $T=(T_1, T_2, ..., T_n)$ is an isomorphism from $(X_1, \mathcal{B}(X_1), \mu_1)$  onto $(X_1 \times Y, \mathcal{B}(X_1 \times Y), \gamma)$, it follows that

\begin{equation*}
 \int_{X_1} c(T_1 x_1, T_2 x_1,..., T_n x_1) \, d\mu_1=\sum_{i=1}^n\int_{X_1} \phi_i(T_ix_1) \, d \mu_1.
 \end{equation*}
 from which together with the fact that $\sum_{i=1}^n \phi_i(x_i)  \leq c(x_1,..., x_n)$ we obtain
 \begin{equation} \label{fc0}
 c(T_1 x_1, T_2 x_1,..., T_n x_1)=\sum_{i=1}^n \phi_i(T_ix_1).  \qquad \mu_1-a.e.
 \end{equation}
 Since $\phi_1$ is $\mu$ almost surely differentiable and ${T_1}_{\#}\mu_1=\mu_1, $ it follows that
 \begin{equation} \label{fccc}
D_1 c(T_1 x_1, T_2 x_1,..., T_n x_1)=\nabla \phi_1(T_1x_1)  \qquad \mu_1-a.e.
 \end{equation}
where $D_1 c $ stands for the partial derivative of $c$ with respect to the first variable.
 Let $A_\gamma \in \mathcal{B}(X_1)$ be the  set with $\mu(A_\gamma )=1$ such that (\ref{fc0}) and (\ref{fccc}) hold for all $x_1 \in A_\gamma, $ i.e.
 \begin{equation} \label{fc}
c(T_1 x_1, T_2 x_1,..., T_n x_1)=\sum_{i=1}^n \phi_i(T_ix_1) \quad \& \quad D_1 c(T_1 x_1, T_2 x_1,..., T_n x_1)=\nabla \phi_1(T_1x_1)  \qquad \forall x_1 \in A_\gamma.
 \end{equation}
Since $\mu_1(X_1\setminus A_\gamma)=0,$ it follows from (\ref{kre}) that \[\int_{ext \, M(T_1, \mu_1)} \varrho(X_1\setminus A_\gamma)\, d\xi (\varrho)=\mu(X_1\setminus A_\gamma)=0,\]
and therefore
there exists  a $\xi$-full measure subset $K_\gamma$ of $ext \, M(T_1, \mu_1)$ such that $\varrho(X_1\setminus A_\gamma)=0$ for all $\varrho \in K_\gamma.$
Let us now define
\[\mathcal{K}:=\big \{F\in \mathcal{S}(T_1);\,\, \exists \varrho \in K_\gamma \text{ with } \mu=F_\# \varrho   \big \},\]
where $\mathcal{S}(T_1)$ is the set of all sections of $T_1.$ Let $\mathcal{G}(\mathcal{K})$
be the set of all sections of $T_1$ generated by $\mathcal{K}$ as in Definition \ref{gen}.  By Proposition \ref{vint}, 
 there exists a  sequence $\{F_i\}_{i=1}^k \subset \mathcal{G}(\mathcal{K})$ with  $ k \in \mathbb{N} \cup\{\infty\}$ satisfying conditions $(i), (ii)$ and $(iii)$ in that proposition. \\

{\it Claim.} We  have  that $k \leq m.$ \\
To prove the claim assume that $k>m$ and 
for each $1\leq n\leq m+1$ let  $\varrho_n={F_n}_\# \mu.$
 It follows  from (\ref{fc}) that
 \begin{equation}
D_1 c\big (T_1\circ F_n (x_1), T_2 \circ F_n(x_1),...,T_n \circ F_n(x_1) \big)=\nabla \phi \big(T_1\circ F_n(x_1)\big)  \qquad  \forall x_1 \in F_n^{-1}(A_\gamma).
 \end{equation}
It follows that
\begin{equation}\label{fc3}
D_1 c\big (x_1, T_2 \circ F_n(x_1),...,T_n \circ F_n(x_1) \big)=\nabla \phi \big(x_1)  \qquad  \,\,\forall x_1 \in \cap_{n=1}^{m+1}F_n^{-1}(A_\gamma),  \,\,  \forall n\leq m+1.
 \end{equation}
 Since $\varrho_n(X_1 \setminus A_\gamma)=0$ and $\varrho_n$ is a probability measure we have  that $\varrho_n(A_\gamma)=1.$
Note that $\varrho_n(A_\gamma)=\mu\big(F_n^{-1}(A_\gamma)\big)$ and therefore $\mu_1 \big(\cap_{n=1}^{m+1}F_n^{-1}(A_\gamma)\big)=1.$   This together with (\ref{fc3}) yield that

  \begin{equation}\label{ku}
D_1 c\big (x_1, T_2 \circ F_n(x_1),...,T_n \circ F_n(x_1) \big)=\nabla \phi \big(x_1) \qquad  \,\,\forall x \in \bar A_\gamma,
 \end{equation}
 where $\bar A_\gamma=\cap_{n=1}^{m+1}F_n^{-1}(A_\gamma).$
 Note that by condition $(i)$ in Proposition \ref{vint} we have  $\mu(B_{m+1})>0.$
 Take $x_1 \in \bar A_\gamma\cap B_{m+1}.$
 It follows from  the $m$-twist condition on splitting sets that the cardinality of 
 \[L_{x_1}:=\Big \{(x_1,y) \in S; \, D_1c\big (x_1, T_2 \circ  F_1(x_1),..., T_n \circ F_n(x_1)\big )=D_1c\big (x_1, y\big ) \Big\},\]
 is at most $m$. On the other hand it follows from (\ref{ku}) that $\big(x_1,  T_2 \circ F_n(x_1),..., T_n \circ F_n(x_1)\big) \in L_{x_1}$ for all $n \in \{1,...,m+1\}.$
 Thus, there exist $i< j\leq m+1$ such that $\big( T_2 \circ F_i(x_1),..., T_n \circ F_i(x_1)\big) =\big( T_2 \circ F_j(x_1),..., T_n \circ F_j(x_1)\big) .$
 Since  $T_1 \circ  F_i=T_1 \circ F_j=Id_{X_1}$ and  the map $T=(T_1,T_2,..., T_n)$ is injective it follows that $F_{i}(x_1)=F_j(x_1).$ 
This is a contradiction as $x_1 \in B_j \subseteq B_{m+1}$ from which the claim follows.\\
  By the latter  claim  we have that $k\leq m.$ It then follows from Proposition \ref{vint} that  every 
 $\varrho \in K_\gamma$ is absolutely continuous with respect to the measure $\sum_{i=1}^k \varrho_i  $ where $\varrho_i={F_i}_\# \mu$ for $1 \leq i \leq k.$
 For every $B \in \mathcal{B}(X)$ it follows from (\ref{kre}) that
 \begin{eqnarray*}\mu_1(B)=\int_{ext \, M(T_1, \mu_1)}  \varrho(B)\, d\xi (\varrho)= \int_{K_\gamma} \varrho(B)\, d\xi (\varrho),
\end{eqnarray*}
from which we obtain that $\mu$ is absolutely continuous with respect to $\sum_{i=1}^k\varrho_i.$ It follows  that $d\mu/d(\sum_{i=1}^k\varrho_i)=\alpha(x)$ for some measureble non-negative
 function $\alpha.$ 
 Assume that  $F_1,..., F_{k}$ are $(\mathcal{B}(X_1)_\mu,\mathcal{B}(X))$-measurable sections  of
 the mapping $T_1$ with $\varrho_i={F_i}_\# \mu_1.$  Setting $\alpha_i=\alpha \circ F_i,$ it follows from  ${T_1}_\# \mu_1=\mu_1$ that $\sum_{i=1}^k \alpha_i(x)=1$ for $\mu_1$-a.e. 
$x \in X_1.$  It also follows from Corollary 6.7.6 in \cite{Bo} that each $F_i$ is $\mu_1$-a.e. equal to  a 
$(\mathcal{B}(X_1),\mathcal{B}(X_1))$-measurable
 function  still denoted  by $F_i.$  For each $i \in \{1,...,k\},$ let  $G_i=\big (T_2 \circ F_i,..., T_n\circ F_i\big).$  We now show that
 $\gamma =\sum_{i=1}^k \alpha_i (\text{Id}\times  G_i)_\# \mu$.  For each bounded continuous function $f : X_1 \times Y \to \R$ it follows that
\begin{eqnarray*}
\int_{X_1 \times Y} f(x,y) \, d \gamma=\int_{X_1} f(T_1x_1,T_2x_1,...,T_nx_1) \, d\mu_1&=&\sum_{i=1}^k\int_{X_1} \alpha(x_1) f(T_1x_1,T_2x_1,...,T_nx_1)  \, d\varrho_i\\
&=& \sum_{i=1}^k\int_{X_1} \alpha\big(F_i (x_1) \big)  f\big (T_1\circ F_i (x_1),T_2 \circ F_i (x_1),...,T_n \circ F_i (x_1)  \big ) \, d\mu_1\\
&=& \sum_{i=1}^k\int_{X_1}\alpha_i(x) f\big (x_1,G_i (x_1) \big ) \, d\mu_1.
\end{eqnarray*}
Therefore,
$\gamma =\sum_{i=1}^k \alpha_i (\text{Id}\times   G_i)_\# \mu.$ This completes the proof.
\hfill $\square$\\

Proof of Proposition \ref{mm}
 goes in the same lines as the proof of Theorem \ref{main} and the only difference is that  by using the same argument one obtains $k\in \mathbb{N} \cup\{\infty\}$ instead of being bounded by $m$ as in the case of the $m$-twist condition. \\

We conclude this section by proving the generalized-twist property for the  locally 1-twisted costs.\\

\textbf{Proof of Proposition \ref{non}}. Assume that $S \subset X_1 \times ...\times X_n$ is a $c$-splitting set.  Fix $(\bar x_1, ...,\bar x_n) \in S$. We need to show that
the set \[L=\Big\{(\bar x_1, x_2,...,x_n) \in S; \, D_1 c(\bar x_1 , \bar x_2,...,\bar x_n)=D_1 c(\bar x_1, x_2,...x_n)\Big\},\]
is finite. If $L$ is not finite there exists an infinitely countable subset $\{(\bar x_1, x^k_2,...x^k_n)\}_{k \in \mathbb{N}} \subset L.$  
Since $S$ is  compact  then the sequence  $\{(\bar x_1, x^k_2,...x^k_n)\}_{k \in \mathbb{N}}$ has an
accumulation point $(\bar x_1, x^0_2,...x^0_n)\in S$ and there exists a subsequence still denoted by $\{(\bar x_1, x^k_2,...x^k_n)\}_{k \in \mathbb{N}}$ 
such that $x_i^k \to x_i^0$ as $k \to \infty$ for $i=2,...,n.$  Since
 $D_1 c$ is continuous it follows that
$(\bar x_1, x^0_2,...x^0_n) \in L.$ Since $c$ is  locally 1-twisted on $S$, this leads to a contradiction   as $(\bar x_1, x^0_2,...x^0_n) $ is an accumulation point of the
sequence $\{(\bar x_1, x^k_2,...x^k_n)\}_{k \in \mathbb{N}}$ and 
\[D_1 c(\bar x_1, x^0_2,...x^0_n)=D_1 c(\bar x_1, x^k_2,...x^k_n), \qquad \forall  k \in \mathbb{N}.\] This completes the proof. \hfill $\square$

\end{document}